\documentclass[runningheads]{llncs}
\usepackage{graphicx}
\usepackage{amsmath, amssymb}
\usepackage{float}

\usepackage{booktabs}  
\usepackage{color}

\newcommand{\black}[1]{\textcolor{black}{#1}}

\begin{document}
\title{A Computational \black{Comparison} of  Optimization Methods for the Golomb Ruler Problem}
\titlerunning{A Computational \black{Comparison} of  Optimization Methods  for Golomb Ruler}
%
\author{Burak Kocuk\inst{1}\orcidID{0000-0002-4218-1116} \and
Willem-Jan van Hoeve\inst{2}\orcidID{0000-0002-0023-753X}}
\authorrunning{Kocuk and van Hoeve}
%
\institute{Sabanc{\i}  University, Istanbul, 34956 Turkey \and
Carnegie Mellon University, Pittsburgh, PA, 15213 USA\\
\email{burakkocuk@sabanciuniv.edu, vanhoeve@andrew.cmu.edu}}
\maketitle              
\begin{abstract}
The Golomb ruler problem is defined as follows: Given a positive integer $n$, locate $n$ marks on a ruler such that  the distance between any two distinct pair of marks are different from each other and the total length of the ruler
is minimized. The Golomb ruler problem has applications in information theory, astronomy and communications, and it can be seen as a challenge for combinatorial optimization algorithms. Although constructing high quality  rulers is well-studied, proving optimality is a far more challenging task.
In this paper, we provide a computational \black{comparison} of different optimization paradigms, each using a different model (linear integer, constraint programming and quadratic integer) to certify that a given  Golomb ruler is optimal.   We  propose several enhancements to improve the computational performance of each method by exploring bound tightening, valid inequalities, cutting planes and branching strategies. We conclude that a certain quadratic integer programming model solved through a Benders decomposition and strengthened by two types of valid inequalities performs the best in terms of solution time for small-sized Golomb ruler problem instances. On the other hand, a constraint programming model improved by   range reduction and a particular branching strategy could have more potential to solve larger size instances due to its promising parallelization features.

\end{abstract}

\section{Introduction}
For a given positive integer $n$, let us denote the positions of $n$ marks on a ruler as $x_1, x_2, \dots, x_n$. Without loss of generality, we assume that the position of the first mark is zero, i.e. $x_1=0$, and the locations are ordered, i.e. $x_1 \le x_2 \le \dots \le x_n$. A \textit{Golomb ruler} satisfies the property that the pairwise distances between distinct marks are all different, in other words, $x_j-x_i \neq x_k-x_l$ unless $i=l$ and $j = k$. The \textit{optimal} Golomb ruler is the one with the smallest length, that is, a Golomb ruler with the minimum $x_n$.
The Golomb ruler problem  has interesting applications in several fields \cite{Bloom1977}, including information theory \cite{Robinson1967},
astronomy, 
and communications \cite{babcock1953,Oshiga2014,Blum1974,Blum1975}. 

In general, constructing a Golomb ruler with a given number of marks
is an easy task, and there are many heuristic methods that provide
high quality rulers. For instance, previous literature on
heuristics has focused on affine and projective plane constructions
\cite{singer1938,drakakis2006}, genetic {algorithms}
\cite{soliday1995}, and local search \cite{prestwich2001,dotu2005}
while exact methods based on constraint programming
\cite{smith1999,galinier2003} {or hybrid methods \cite{slusky2013}}
exist as well. Although not proven to be {NP-hard yet}, solving the
Golomb ruler problem exactly proved to be notoriously difficult. For
instance, the optimal rulers for $n=24,25,26,27$ have been proven by a
parallel search with thousands of workstations coordinated by the
website \texttt{distributed.net}, and it took approximately 4, 4, 1,
and 5 years to complete, respectively. Currently, a search for the
28-mark problem is under way for more than 4 years.

As summarized above, most of the effort to prove that a given Golomb ruler is an optimal one is devoted to explicit enumeration techniques. 
However, such brute force approaches seem to be the only viable option since it is very difficult to establish strong valid lower bounds for the Golomb ruler problem.

At this point, we would like to state the main purpose of this paper, which is to \black{\textit{certify the optimality of a given Golomb ruler through optimization methods}}. Most optimization algorithms inherently solve relaxations and hence, naturally provide lower bounds for minimization type problems. Therefore, it is worth focusing on  optimization models to better understand the structure of the Golomb ruler problem, and hopefully, {propose efficient methods} which we can use to solve the Golomb ruler problem instances.

In this paper, we consider three classes of optimization problems to
carry out the aforementioned analysis. In particular, we formulate the
Golomb ruler problem as linear integer programming, constraint
programming and quadratic integer programming problems. Some of these
models exist in the literature {while and the others are
  introduced}, to the best of our knowledge, by us. Since the
performance of the basic models is not satisfactory to solve instances
with more than 10 marks, we propose several enhancements to
improve the scalability of each method {by means} of bound
tightening, valid inequalities, cutting planes and effective search
strategies. Our computational experiments show that linear
  integer programming models scale up to {13 marks} given a budget of 8
  hours while constraint programming {models} can solve up to 13-mark
  instances in about an hour and 14-mark problem in about 10
  hours. Quadratic integer programming models, on the other hand, are
  able to solve 14-mark instance in about four hours,
 all using a modest personal computer.  
 \black{As a comparison,  it took {2.8} hours for the constraint programming model in~{\cite{galinier2003}} to find an optimal ruler for the 13-mark instance and another {11.8} hours to prove its optimality.
The lean implementation of the search method  in~{\cite{slusky2013}} reduced the respective computational effort to 0.6 and 1.3 hours, albeit at the expense of a significantly larger search tree.}


The rest of the paper is organized in three sections, which
respectively cover {the} linear integer programming, constraint
programming and quadratic integer programming models in detail. Each
section contains a formulation, enhancements and computational
experiments subsections together with some discussions and comparisons
of these different optimization paradigms. 

\section{Linear Integer Programming Models}

\subsection{Two  Formulations}
In this section, we present two {known} linear integer programming
formulations for the Golomb ruler problem with $n$ marks. One of these
formulations is \textit{exact} and the other one is only a
\textit{relaxation} but can be made exact by the use of additional
features as detailed later. For both of the models, we assume that an
upper bound, say $L$, on the optimal length is known (such a bound can be obtained as the length of any feasible Golomb ruler with $n$ marks).

\subsubsection{``$d+e$" Formulation}
We first present an exact linear integer programming model for the Golomb ruler problem~{\cite{lorentzen1991}}. In this formulation, there are two sets of decision variables: Let $e_{ijv}$ be one if the distance between marks $i$ and $j$ is $v$, and zero otherwise. Also, we define $d_{ij}$ as the distance between marks $i$ and $j$. Then, the optimization problem is given as follows:
\begin{subequations} \label{eq:lp1}
\begin{align}
 \min_{d, e}  & \ \sum_{i=1}^{n-1}   d_{i,i+1} \label{eq:ilp:obj} \\
\text{s.t.} & \ \sum_{v=1}^L e_{ijv} = 1  & 1& \le i < j \le n  \label{eq:ilp:assign} \\
& \ \sum_{i<j} e_{ijv} \le  1 & v&=1,\dots, L  \label{eq:ilp:golomb}  \\
& \ \sum_{v=1}^L v e_{ijv} = d_{ij} & 1& \le i < j \le n  \label{eq:ilp:definition} \\
& \ \sum_{k=i}^ {j-1} d_{k,k+1} = d_{ij} & 2 &\le i + 1 < j \le n   \label{eq:ilp:identity} \\
& \  e_{ijv} \in \{ 0,1\}, d_{ij} \in \mathbb{Z}_+ & 1& \le i < j \le n,   v=1,\dots, L.
\end{align}
\end{subequations}
Here, the objective function~\eqref{eq:ilp:obj} minimizes the length of the Golomb ruler (alternatively, it can be simply given as $d_{1,n}$). Constraint ~\eqref{eq:ilp:assign} assigns a distance between $1, \dots, L$ to every pair of marks  $i$ and $j$ while constraint  ~\eqref{eq:ilp:golomb} ensures that each distance between $1, \dots, L$ is assigned at most once.  Constraint~\eqref{eq:ilp:definition} is simply a definition of the $d$ variables in terms of the $e$ variables. Finally,  constraint~\eqref{eq:ilp:identity} is an identity guaranteeing that the distance between the marks $i$ and $j$ is the sum of the basic distances between consecutive marks. We note that $d$ variables can be projected out by substituting the definition given in~\eqref{eq:ilp:definition} into constraint~\eqref{eq:ilp:identity} although the resulting lower dimensional model does  not seem to be more advantageous empirically in terms of computational efficiency.

\subsubsection{``$d$" Formulation}
Now, we present a relaxation of the Golomb ruler problem by
eliminating the $e$ variables from the ``$d+e$" Formulation. The
resulting optimization problem is defined as follows:
\vspace{-2mm}
\begin{subequations}\label{eq:relax}
\begin{align} 
\min_d  & \ \sum_{i=1}^{n-1}   d_{i,i+1} \\
\text{s.t.}  & \ \eqref{eq:ilp:identity} \notag \\
 & \ \sum_{(i,j)\in R} d_{ij} \ge \frac12|R|(|R|+1) &R& \subseteq \{(k,l)\in\mathbb{Z}^2: 1\le k<l \le n\} \label{eq:sumDist}.
\end{align}
\end{subequations}
{The constraint set~\eqref{eq:sumDist}, called ``Subset Sum
Inequalities'', was introduced in \cite{lorentzen1991} to strengthen
the integer programming (IP) formulation, and} defines the facets of the convex hull of the
all-different constraint \cite{williams2001representations}. In
practice, the model~\eqref{eq:relax} can be solved with a constraint
generation scheme in which the subset sum inequalities are gradually
added. We also note that the complexity of the separation of these
inequalities is polynomial-time as it requires sorting ${n \choose 2}$
many numbers.

As opposed to the ``$d+e$" Formulation,  the ``$d$" Formulation is not an exact representation of the Golomb ruler problem.
However, it is known that the optimal value of the ``$d$" Formulation is equal to the linear programming (LP) relaxation value of ``$d+e$"  Formulation \cite{meyer2006}.

We note that even if the ``$d$" Formulation is solved as an IP, it is still not an exact formulation since it does not guarantee that $d_{ij} \neq d_{kl}$ for $i\neq k$ and $j \neq l$, i.e., the uniqueness of the distances.
However, this observation leads to a natural way  to make  the ``$d$" Formulation exact: We can  solve the problem~\eqref{eq:relax} as an LP with two callbacks: Firstly, we add  lazy constraint callbacks to ensure that subset sum inequalities~\eqref{eq:sumDist} are satisfied. Secondly, we add a branch callback such that the missing constraint $d_{ij} \neq d_{kl}$ is enforced by the solver as we go down the branch-and-bound tree. In particular, once  $d_{ij} = d_{kl}$ for $i\neq k$ and $j \neq l$, we can create a dichotomy as $d_{ij} \le d_{kl} - 1$  and  $d_{ij} \ge d_{kl} + 1$. 

\subsection{Enhancements}
We propose several enhancements to the models presented above based on bound tightening and branching strategies. We also carried out a preliminary polyhedral study of the Golomb ruler problem with the hope of obtaining strong valid inequalities. 
{Although we have discovered several families of valid inequalities, they have not helped solving the problem more efficiently. Therefore, we leave this line of research as future work for further inquiry.
}

In the sequel, let $G_m$ denote the length of the optimal Golomb ruler of order $m$, $m=1,\dots,n-1$. We assume that all the $G_m$ values are known for $m < n$ when we are trying to solve the $n$-mark problem.

\subsubsection{Bound Tightening}
\label{sec:bound lp}
Bound tightening is  a widely used strategy in optimization algorithms to reduce the range of the decision variables in order to save computational effort. It can also be used as a way to strengthen relaxations and improve the performance of search methods \black{(see \cite{galinier2003} for an application  to the Golomb ruler problem}). Our bound tightening procedure starts with the following simple observation: If the difference between marks $i$ and $j$  is small (large), then $d_{ij}$ cannot be too large (small). In particular, we can infer the following initial bounds on $d_{ij}$ variables:
\begin{equation}
G_{j-i+1} =: \underline d_{ij} \le d_{ij} \le  \overline d_{ij} :=  L - G_{i} - G_{n-j+1}. \label{eq:golombbounds}
\end{equation}

After  this initialization step, we can further improve the bounds $\underline d_{ij}$ and $\overline d_{ij}$ by solving bounding LPs. In particular, we can  minimize/maximize $d_{ij}$ variables over a suitable relaxation (for instance, over the feasible region of the LP relaxation of~\eqref{eq:lp1}) to try  to improve the bounds iteratively. 
Once new bounds $\underline d_{ij} $ and $ \overline d_{ij} $ are obtained after rounding up and down the minimum and maximum values, we  repeat this procedure until the fixed point is reached, that is, none of the bounds are improved further.  We will refer to this procedure as \textit{LP Bounding}.

As an additional mechanism to tighten the variable bounds on the $d$ variables, we extend the LP Bounding approach in the following sense. Now, we optimize the $d$ variables over the feasible region  of~\eqref{eq:lp1} (not its LP relaxation as in the LP Bounding approach) with a limited computational budget. We will refer to this procedure as \textit{IP Bounding}. Although this approach requires an additional non-trivial effort, it pays off in terms of reducing the solution time of the ``$d+e$'' formulation.

We  also use the bounds $\underline d_{ij}$ and $\overline d_{ij}$ to fix some of the binary variables $e_{ijv}$ to zero as follows:
\begin{equation*}
e_{ijv} = 0 \text{ if } v < \underline d_{ij} \text{ or } v > \overline d_{ij}.
\end{equation*}
This procedure reduces the total range of the $d$ variables and the number of $e$ variables considerably although it does not improve the LP relaxation bound.

\subsubsection{Branching Strategies}
\label{sec:lip branch}

The choice of branching strategies may significantly affect the computational performance of the mixed-integer programming solvers. Branching decisions made by the solvers can be altered by either explicitly choosing the variables to be branched on through branch callbacks, or implicitly by assigning priorities to the integer variables. We experimented with both of these choices by exploiting the  structure of the Golomb ruler problem.

In terms of imposing explicit branching decisions, we experimented with two strategies which can be applied to both  ``$d+e$" and    ``$d$" Formulations. The first strategy, which we will call as ``Left Branching", is described as follows: We first branch on the variable $d_{12}$ by creating $\overline d_{12} - \underline d_{12} + 1$ many child nodes, each taking an integer from the interval $[\underline d_{ij}, \overline d_{ij} ]$. 
Then, we proceed by solving the node relaxations. Whenever we have to make another branching decision, we decide on the next variable still undecided from the left of the ruler (for instance, the second variable would be $d_{23}$). Such an algorithm is based on the intuition that the classical dichotomy branching is probably ineffective for the Golomb ruler problem since assigning particular values to $d_{ij}$ variables can help detecting infeasibility of more branches than simply branching on $d_{ij}$ by assigning intervals, or traditional binary branching on the $e_{iju}$ variables. The shortcoming of this strategy is of course the increased number of child nodes for each level of the branch-and-bound tree.

The second strategy for imposing explicit branching decisions is called ``Difference Branching", 
and implemented with the inclusion of a branch callback function. 
  Suppose that in a certain node in the branch-and-bound tree, we have two variables $d_{ij}$ and $d_{kl}$, $i\neq k$ and $j \neq l$, such that $|d_{ij} - d_{kl}| \le 1$. Then, we can branch on constraints $d_{ij} \le d_{kl} - 1$  and  $d_{ij} \ge d_{kl} + 1$ as this is a valid partitioning of the feasible region.

In terms of imposing implicit branching decisions, we experimented with different priority assignment strategies. The  most successful strategy seems to  be the one that assigns higher priorities to $e_{iju}$ variables with smaller $u$ indices. Here, the intuition is that if smaller distances are decided first,  then we can either find feasible solutions or detect infeasibility faster.


\subsection{Computations}
\label{sec:d+e comp}

%
%

We first report the  results of our computational experiments with the ``$d+e$'' formulation in Table~\ref{tab:d+e} for $n=9,\dots,14$.  We compare the following two settings \black{(we note that the bound tightening techniques are only applied at the root node assuming that  $L=G_n$)}:
\begin{itemize}
\item
LP Bounding: Bound tightening is applied over the LP Relaxation of the   ``$d+e$'' formulation in parallel  for five rounds.
\item
IP Bounding after LP Bounding: Bound tightening is applied over the    ``$d+e$'' formulation with a budget of 1 second for $n\le 12$ and 1 minute for $n \ge 13$ in parallel  for five rounds.
\end{itemize}
\begin{table}[t]\footnotesize
\centering
\caption{Results of the ``$d+e$'' Formulation with LP  and IP Bounding strategies. Here, BTT, OT, TT and \#BBN respectively represent  the bound tightening, optimization, total time (in seconds unless otherwise stated), and the number of branch-and-bound nodes.}\label{tab:d+e}
\begin{tabular}{rcrrrrcrrrr}
\toprule
           &&                    \multicolumn{ 4}{c}{LP Bounding} &&                 \multicolumn{ 4}{c}{IP Bounding after LP Bounding} \\
\cmidrule{3-6} \cmidrule{8-11}
        $ n$ &~~~&        BTT &         OT &         TT &    \#BBN &~~~~&        BTT &         OT &         TT &   \#BBN \\
\midrule
%


         9 &&       0.10 &       0.36 &       0.46 &          0 &&      41.35 &       0.60 &      42.05 &        615 \\

        10 &&       0.11 &       1.42 &       1.53 &        333 &&      37.69 &       0.68 &      38.48 &        337 \\

        11 &&       0.12 &     199.87 &     199.99 &     155,421 &&      20.55 &     140.38 &     161.05 &     138,820 \\

        12 &&       0.19 &     297.62 &     297.81 &     153,244 &&      68.27 &     208.39 &     276.85 &     126,405 \\

        13 &&       0.26 &   32,435.04 &   32,435.30 &   13,949,679 &&    2,874.62 &   24,993.13 &   27,868.01 &   13,363,776 \\

        14 &&       0.32 &      $ >10$ hr&          - &          - &&    3,721.56 &            $ >10$ hr &          - &          - \\
\bottomrule
\end{tabular}  
\end{table}
CPLEX 12.8 is used as the mixed-integer linear programming (MILP) solver on a 64-bit computer with Intel Core i7
CPU 2.60GHz processor and 16 GB RAM. \black{Since our aim to prove the optimality of the $n$-mark ruler of length $G_n$, we set $L=G_n-1$ while solving the MILPs.}

We   note that the LP Bounding scheme is quite cheap to obtain better variable bounds than the initial bounds derived in~\eqref{eq:golombbounds}. On the other hand,  IP Bounding requires an additional nontrivial  effort to further improve those bounds. We observe that this additional computational effort can be justified when $n \ge 12$ as the reduction in the optimization step overweighs the increase in the bound tightening step and the 13-mark instance can be solved in less than 8 hours in total. Due to the sharp increase in the CPU time necessary, we were not able to solve the Golomb ruler problem with 14 marks in less than 10 hours.

{We report} the results of our experiments with the ``$d$''
formulation in Table~\ref{tab:d} for $n=9,10,11$. We compare the
Difference (Diff.) Branching and Left Branching strategies as introduced in
Section \ref{sec:lip branch} in combination with LP Bounding and IP
Bounding. We observe that the Left Branching becomes significantly
better than the Difference Branching approach as the number of marks
increases. Since we were able to solve only very small Golomb ruler
problem instances with the ``$d$'' Formulation in comparison to the
``$d+e$'' Formulation, we have not pursued this direction further.
Nevertheless, the Left Branching strategy has proved to be quite
effective and is utilized multiple times in {this} paper.
\begin{table}[t]\footnotesize
\centering
\caption{Results of the ``$d$'' Formulation with different variable bounding and branching strategies.} \label{tab:d}
\begin{tabular}{rcrrcrrcrrcrr}
\toprule
           &~&                 \multicolumn{ 5}{c}{LP Bounding} &~~& \multicolumn{ 5}{c}{IP Bounding after LP Bounding} \\
\cmidrule{3-7} \cmidrule{9-13}
           && \multicolumn{ 2}{c}{Diff. Branching} &~~ & \multicolumn{ 2}{c}{Left Branching} & & \multicolumn{ 2}{c}{Diff. Branching} &~~&  \multicolumn{ 2}{c}{Left Branching} \\
\cmidrule{3-4} \cmidrule{6-7} \cmidrule{9-10} \cmidrule{12-13} 
        $ n$&&         TT &    \#BBN &&         TT &    \#BBN &&         TT &    \#BBN &&         TT &    \#BBN \\
\midrule

         9 &&       1.45 &       3,277 &&       3.05 &       5310 &&      42.59 &       3,257 &&      44.45 &       5,338 \\

        10 &&       0.55 &        805 &&       2.00 &       3248 &&      38.21 &        785 &&      39.74 &       3,297 \\

        11 &&    9,915.67 &    4,256,165 &&    1,920.96 &    1,651,695 &&    9,785.72 &    4,285,299 &&    1,862.69 &    1,655,038 \\
        
\bottomrule
\end{tabular}  
\end{table}

\section{Constraint Programming Model}
In the previous section, we presented {two} linear integer
programming models to solve the Golomb ruler problem. Although several
enhancement of these models are introduced and the computational
effort is reduced significantly, we were not able prove the optimality
of a given 14-mark ruler in a reasonable amount of time. In this
section, we switch our attention to constraint programming models
which prove to be more successful for the Golomb ruler problem.

\subsection{Formulation}
\label{sec:cp form}
A constraint programming model of the Golomb ruler problem can be
formulated as {follows~\cite{smith2000,slusky2013}}:
\begin{subequations} \label{eq:cp}
\begin{align}
 \min_d  & \  d_{1,n} \\
\text{s.t.} & \ \text{alldiff}(\{d_{ij}: \  1 \le i < j \le n\})   \label{eq:cp alldiff} \\
& \ d_{ij} + d_{jk} = d_{ik} & 1& \le i < j < k \le n \label{eq:cp triangle} \\
& \  d_{ij } \in \{ \underline{d}_{ij},  \dots, \overline{d}_{ij}\} & 1& \le i < j \le n.  \label{eq:cp bounds}
\end{align}
\end{subequations}
Here, constraint~\eqref{eq:cp alldiff} ensures that each distance $d_{ij}$ are different from each other. Constraint~\eqref{eq:cp triangle} guarantees that the distances respect the ``triangle" constraint, that is, the distance between marks $i,k$ is the sum of the distances between marks $i,j$ and $j,k$, where $j$ is strictly between $i$ and $k$. Finally, constraint ~\eqref{eq:cp bounds} specifies the ranges of the decision variables.

\subsection{Enhancements}
\label{sec:cp enh}

The constraint programming model~\eqref{eq:cp} can easily solve small instances of the Golomb ruler problem, e.g., $n \le 10$, but runs into slow convergence {issues} for even slightly larger instances.  Similar to the integer programming models considered in the previous sections, we propose some enhancements to improve the scalability of the constraint programming model. These enhancements utilize  bound tightening,  table constraints and search strategies.

\subsubsection{Bound Tightening}

The  bound tightening procedures proposed in Section \ref{sec:bound lp} based on LPs and IPs are quite effective in reducing the ranges of the $d$ variables. We now discuss another similar procedure based on constraint programming techniques \black{(see \cite{martin1996new} for a related method called ``shaving'')}. The proposed idea is quite simple: We fix a $d_{ij}$ variable to its current lower or upper bound and solve the feasibility version of the constraint programming model~\eqref{eq:cp} for a limited amount of time. If the infeasibility of this restricted model can be proven, this implies that we can tighten the range of the $d_{ij}$ variable by excluding the value that we have fixed. We will refer to this procedure as \textit{CP Bounding}. This approach is implemented in an iterative fashion with limited computational budget and proved to be helpful to further reduce the range of the decision variables.

\subsubsection{Table Constraints}

Table constraints can be crucial to speed up constraint programming
solvers. By exploiting the specific structure of the Golomb ruler
problem, we also define certain ``forbidden assignments''. The
construction is as follows: Consider the subruler with marks numbered
as {$\{i, \dots, i+4\}$}. We first enumerate the set $S_j$ of all triplets 
{$(d_{j,j+1}, d_{j+1,j+2}, d_{j+2,j+3})$}
that {constitute a feasible subruler with respect} to the variable bounds and
all-different constraints{, for $j=i$ and $j=i+1$}.
Now, suppose that for some
$(\bar d_{i,i+1}, \bar d_{i+1,i+2}, \allowbreak \bar d_{i+2,i+3})\in
S_i$, there does not exists any $\bar d_{i+3,i+4}$ such that
$( \bar d_{i+1,i+2}, \bar d_{i+2,i+3}, \allowbreak \bar d_{i+3,i+4})$
belongs to the set $ S_{i+1}$. In this case, we can declare the
triplet
$(\bar d_{i,i+1}, \allowbreak \bar d_{i+1,i+2}, \bar d_{i+2,i+3})$ as
``forbidden''. We can repeat this procedure for a few rounds across
different subrulers to identify more forbidden assignments.

\subsubsection{Search Strategies}
 
Search strategies are extremely important in constraint programming as they significantly alter the performance of the solvers. Inspired by the Left Branching for the linear integer programming model and its adaptation to the quadratic integer model, we have decided to employ a variable selection rule based on lexicographical ordering. We also set the search phase parameter to depth first search as our aim is to prove the optimality of a given ruler efficiently. Finally, we experiment with different value selection strategies and decide to use  the one based on the smallest impact.

\subsection{Computations}
\label{sec:cp comp}

We report the results of our computational experiments with the
constraint programming formulation in Table~\ref{tab:cp} for
$n=9,\dots,14$. In addition to the ``LP Bounding'' and ``IP Bounding
after LP Bounding'' settings introduced in Section~\ref{sec:d+e comp},
we also experimented with the following version:
\begin{itemize}
\item CP Bounding after LP and IP Bounding: Bound tightening is
  applied over the constraint programming formulation with a budget of
  1 second for $n\le 12$ and 1 minute for $n \ge 13$ in parallel for
  five rounds.  {This includes the generation of forbidden
    assignments based on the table constraints.}
\end{itemize}
CPLEX \black{CP Optimizer} is used as the constraint programming solver \black{with the default settings unless otherwise stated}.

\begin{table}[t]\footnotesize
\centering
\caption{Results with LP Bounding, IP Bounding and CP Bounding strategies.}\label{tab:cp}
\begin{tabular}{rcrrcrrcrrr}
\toprule
           &~~& \multicolumn{ 2}{c}{LP Bounding} &~~& \multicolumn{ 2}{c}{IP Bounding after LP} &~~& \multicolumn{ 3}{c}{CP Bounding after LP and IP } \\
\cmidrule{3-4} \cmidrule{6-7} \cmidrule{9-11}
        $ n $  &&         OT &         TT &&         OT &         TT &&        BTT &         OT &         TT \\
\midrule



         9 &&       2.82 &       2.92 &&       0.27 &      41.71 &&       1.03 &       0.06 &      42.54 \\

        10 &&      11.95 &      12.06 &&       0.28 &      38.08 &&       1.17 &       0.07 &      39.04 \\

        11 &&     140.30 &     140.42 &&      19.63 &      40.30 &&      11.85 &       0.54 &      33.07 \\

        12 &&     554.71 &     554.90 &&      23.39 &      91.85 &&      27.77 &       1.52 &      97.75 \\

        13 &&   11,615.19 &   11,615.45 &&    1,780.54 &    4,655.42 &&    1,665.97 &     209.87 &    4,750.71 \\

        14 &&           $ >10$ hr &          - &&   31,839.14 &   35,561.02 &&    3,872.62 &   29,795.25 &   37,389.74 \\
\bottomrule
\end{tabular}  
\end{table}

We now summarize our observations:
Firstly, a comparison with Tables \ref{tab:d+e} and \ref{tab:cp} indicates that the constraint programming formulation takes less time than the ``$d+e$" formulation under the same version of bounding for $n \ge 11$. This allows us to solve the 14-mark problem with the constraint programming approach with IP Bounding in 10 hours, which was not possible with the ``$d+e$" formulation.
Secondly, the overhead of the CP Bounding approach is quite large, hence, the reduction in the optimization time compared to the IP Bounding may not be fully justified always. However, we believe that further inquiry along this direction should be pursued. 
For instance, we extended the CP Bounding approach for the 14-mark problem by allowing 10 minutes of budget for each subruler length. This increases the total bounding time to 18,251 seconds but reduces the optimization time to 16,851 seconds. Although the total time remains more or less unchanged, this additional experiment shows that a carefully executed bounding mechanism may have a potential to be  efficient  overall.

\section{Quadratic Integer Programming Models}

{So far, we presented classical linear integer and constraint
  programming formulations for the Golomb ruler problem. In this
  section, we focus on a less-explored approach based on
  quadratic integer programming.}

\subsection{Two Formulations}
In this section, we discuss two possible quadratic integer programming formulations for the Golomb ruler problem with $n$ marks, one based on an \textit{optimization} model and the other based on a \textit{feasibility} version. To the best of our knowledge, such formulations have not been {proposed before} in the literature. 

Let us define a single set of binary variables $y_l$, which takes value one if there is a mark at location $l$ and zero otherwise, $l=1,\dots,L$. Here, $L$ is again an upper bound on the length of a  shortest Golomb ruler with $n$ marks.

\subsubsection{Optimality Version}
We first present an alternative integer programming formulation of the Golomb ruler problem  as follows:
\begin{subequations} \label{eq:qip opt}
\begin{align}
\min_y  & \ \max_{l=1,\dots,L}  \   l  \black{\times} y_l \label{eq:qip opt obj} \\
\text{s.t.} & \ \sum_{l=0}^{L-v} y_l y_{l+v} \le 1 & v&=1,\dots,L \label{eq:qip opt cons} \\
& \  y_l \in \{ 0,1\} &l& =1,\dots, L.\label{eq:qip opt binary}
\end{align}
\end{subequations}
Here, the objective~\eqref{eq:qip opt obj}  minimizes the position of the last mark on a ruler of length $L$, which corresponds to the length of an optimal ruler. We note that the objective function can be easily linearized using an auxiliary variable and enforcing $L$ additional constraints. Constraint~\eqref{eq:qip opt cons} guarantees that each distance $v$ is used at most once in a feasible solution. Observe that the model~\eqref{eq:qip opt} can be reformulated as a quadratically {constrained} program, which contains two types of nonconvexities, one due to the bilinear inequalities ~\eqref{eq:qip opt cons}, and another due to the integrality of the $y$ variables.

Convexification techniques can be utilized to solve or approximate the nonconvex problem~\eqref{eq:qip opt}. We note that this formulation admits a straightforward semidefinite programming (SDP) relaxation given as follows:
\begin{subequations} \label{eq:qip opt sdp}
\begin{align}
\min_{y,z,Y}  &  \ z \label{eq:qip opt sdp obj} \\
\text{s.t.} & \  l \black{\times} y_l  \le z & l&=1,\dots,L \label{eq:qip opt sdp obj aux} \\
& \ \sum_{l=0}^{L-v} Y_{l,l+v} \le 1 & v&=1,\dots,L \label{eq:qip opt sdp cons} \\
& \ 0 \le y_l \le 1 & l&=1,\dots,L \label{eq:relaxed binary} \\
& \  \begin{bmatrix} 1 & y^T \\ y & Y \end{bmatrix} \succeq 0  \label{eq:psd cons}.
\end{align}
\end{subequations}
 Unfortunately, the dual bound obtained from solving the SDP relaxation~\eqref{eq:qip opt sdp} is extremely weak (for instance, the bound obtained for the 10-mark instance is only 15.14 while the length of the optimal Golomb ruler is 55). 
 Therefore, we have not pursued this line of research direction further.
%
%
%
%
%
%
%

\subsubsection{Feasibility Version}
\label{sec:feasVers}
Now, we consider a ``complementary" version of the problem defined as follows: Given the length of a ruler $L$, maximize the number of marks that can be located onto such a ruler that satisfies the Golomb ruler requirements. This version of the problem can be formulated as follows:
\begin{subequations}  \label{eq:qip feas}
\begin{align}
n_L := \max_y  & \ \sum_{l=0}^L  \    y_l \\
\text{s.t.} & \ \eqref{eq:qip opt cons}-\eqref{eq:qip opt binary}.
\end{align}
\end{subequations}
Note that the formulation~\eqref{eq:qip feas} can be seen as the feasibility version of the model~\eqref{eq:qip opt} in the following sense: If $n_L = n$ but $n_{L-1} = n-1$, then we can certify that $G_n = L$. Hence, in order to obtain the length of a shortest Golomb ruler with $n$ marks, that is $G_n$, we can first solve problem~\eqref{eq:qip feas} with $L= G_{n-1}$, and then increase the value of $L$ until we can locate all of the $n$ marks. Such a procedure gives an indirect way of solving the Golomb ruler problem.

Problem~\eqref{eq:qip feas}  is again a nonconvex, quadratically constrained integer program. Below, we propose two linearization methods that can be used to solve problem~\eqref{eq:qip feas} via an appropriate branch-and-bound method.

\paragraph{Linearization via SDP} The problem~\eqref{eq:qip feas} can be reformulated as a mixed-integer SDP as follows:
\begin{subequations}  \label{eq:qip feas sdp}
\begin{align}
\max_{y, Y}  & \ \sum_{l=0}^L  \    y_l \\
\text{s.t.} & \ \eqref{eq:qip opt sdp cons}, \eqref{eq:psd cons}, \eqref{eq:qip opt binary} \notag.
\end{align}
\end{subequations}
Since the problems in this class are not typically supported by commercial solvers, we implemented our own branch-and-bound algorithm. In this algorithm, we solve the SDP relaxation of the model~\eqref{eq:qip feas sdp}, which replaces the binary restriction~\eqref{eq:qip opt binary} with its continuous relaxation~\eqref{eq:relaxed binary}, at each node of the tree. Our algorithm decides which  $y$ variables to choose for branching, which is discussed in more detail in Section \ref{sec:y form branch}.

\paragraph{Linearization via LP} The problem~\eqref{eq:qip feas} can be also reformulated as a mixed-integer LP as follows:
\begin{subequations}  \label{eq:qip feas lp}
\begin{align}
\max_{y, Y}  & \ \sum_{l=0}^L  \    y_l \\
\text{s.t.} & \ \eqref{eq:qip opt sdp cons}, \eqref{eq:qip opt binary} \notag \\
& \ y_l + y_k  -1 \le Y_{lk}  &l&,k = 1,\dots,L \label{eq:qip opt mcc1}\\
& \ 0 \le Y_{lk} \le y_l  &l&,k = 1,\dots,L. \label{eq:qip opt mcc2}
\end{align}
\end{subequations}
Here, constraints~\eqref{eq:qip opt mcc1}--\eqref{eq:qip opt mcc2} correspond to the McCormick envelopes for the equation $Y_{lk} = y_l y_k$. In general, solving problem~\eqref{eq:qip feas lp} directly as an MILP is quite expensive, partly due to the fact that its LP relaxation is highly degenerate. Therefore, we adopt a Benders decomposition approach, whose problem specific details are presented in Section~\ref{sec:y benders}.

\subsection{Enhancements}

We again propose some enhancements to speed up the solution procedure of the feasibility version of the quadratic formulation of the Golomb ruler problem. In particular, we develop two types of valid inequalities, Benders decomposition for the linearized model~\eqref{eq:qip feas lp} and branching strategies. Improved variable bounds obtained via the bound tightening  procedure presented in Section \ref{sec:bound lp} are also used whenever applicable.

\subsubsection{Valid Inequalities}
In this section, we present two families of valid inequalities, which we {refer to} as ``Golomb" and ``Clique" inequalities. Below, we present their precise formulations together with the intuition behind them.

\paragraph{Golomb Inequalities} Since the number of marks that can be placed onto any subruler of length $t$ is upper bounded by $n_t$, the following inequalities are valid \black{and are added to the root node relaxation}:
\begin{equation}\label{eq:golomb ineq}
\sum_{j=l}^{l + \min\{G_{i+1},L\} + 1} y_j \le i \quad i=2,\dots,n_L; \   l=0,\dots,L- (\min\{G_{i+1},L\} + 1).
\end{equation}
More inequalities of this kind can be obtained as follows: Instead of summing the consecutive~$y$ variables, we can consider any subset of these variables whose indices are separated by exactly the same integer $c$, $c=2,\dots,\lfloor L/2 \rfloor$, such as $y_j, y_{j+c}, y_{j+2c}, \dots$. 

\paragraph{Clique Inequalities} In order to better explain the construction of the clique inequalities, it is more suitable to present the Golomb ruler problem as a special \textit{maximum cardinality clique problem} defined as follows: Let us consider a complete  graph $G=(V,E)$, where the set of vertices is $V:=\{0, \dots, L\}$. We partition the edge set $E$ into $L$ subsets $E_l$ defined as $E_l:=\{(i,i+l): i =0,\dots, L - l\}$ for $l=1,\dots,L$. Then, in order to solve the Golomb ruler problem, we search for a largest clique $G'=(V',E')$ in this graph such that at most one edge from each subset $E_l$ belongs to  $E'$, that is, $|E_l \cap E'| \le 1$.

Motivated by the above construction, let us introduce the clique inequalities, which  are easily implementable in a cutting plane framework. Consider a fractional solution $\tilde y$, and construct two sets $\mathcal{L}_1 :=\{l: \tilde y_l = 1\}$ and $\mathcal{L}_f :=\{l: \tilde y_l \in (0, 1) \}$. Let us define the distances induced by the solution $\tilde y$ as $\mathcal{D} := \{ |k - l|: k,l \in \mathcal{L}_1 \}$. We will now construct an auxiliary graph $\tilde G = (\tilde V, \tilde E)$, where $\tilde V := \mathcal{L}_f $ and 
$\tilde E := \{ (k,l) \in \tilde V \times \tilde V : |k-l| \in \mathcal{D}  \}$. We also associate node weights $\tilde y_l$ for each $l \in \tilde E$. Then, each maximal clique $\tilde C$ in the graph $\tilde G$ whose weight is more than 1 gives rise to a cutting plane of the following form:
\begin{equation}\label{eq:clique ineq}
\sum_{l \in \tilde C} y_l \le 1.
\end{equation} 
We note that the set of all maximal cliques in a graph can be found by the {Bron-Kerbosch algorithm~\cite{bron1973algorithm} in reasonable time for such small graphs,} and the inequalities~\eqref{eq:clique ineq} can be added as local user cuts in a branch-and-cut algorithm.

\subsubsection{Benders Decomposition}
\label{sec:y benders}

Since we observe that solving   problem~\eqref{eq:qip feas lp} directly is not computationally efficient, we employ a Benders decomposition technique instead. In this approach, we solve the following master problem
\begin{subequations}\label{eq:benders master}
\begin{align}
\max_{y}  & \ \sum_{l=0}^L  \    y_l \\
\text{s.t.} & \ \eqref{eq:qip opt binary} \notag \\
 & \ {\sum_{l \in C} |C_u(l)| y_l \le  \frac{|C|}{2} + 1  \qquad C \subseteq \{0,\dots,L\}, u=1,\dots,L}, \label{eq:BendersCutsAll} 
\end{align}
\end{subequations}
where constraints~\eqref{eq:BendersCutsAll} are added in lazy fashion until feasibility is proven. Here, $C_u(l) := \{k: |k-l|=u\} $.
This is achieved through the separation procedure (feasibility check) described as follows: 
Given a binary vector $\tilde y$, we first define the set of marks as
$M:=\{l: \tilde y_l = 1\}$. {For this candidate ruler to be a
Golomb ruler}, the distances between each pair of mark should be
distinct. Therefore, the cardinality of the set
{$M_u(l) := \{k \in M: |k-l|=u \}$} should be~$1$ for $l\in M$ and
$u=1,\dots,L$. Otherwise, we detect infeasibility and can add the
following cut:
\begin{equation}\label{eq:benders cut}
{\sum_{l \in M} |M_u(l)| y_l \le  \frac{|M|}{2} + 1}.
\end{equation}
In other words, we expand the constraint set $\eqref{eq:BendersCutsAll}$ by the inclusion of the set~$M_u$. We note that our approach allows to add multiple cuts of the form~\eqref{eq:benders cut} corresponding to different $u$ values for a given solution.

What we described so far amounts to a classical implementation of the Benders decomposition technique in which a ``multi-tree"  approach is employed, that is, at each iteration, we solve the master problem~\eqref{eq:benders master} as an MILP. Therefore, multiple branch-and-bound trees are created. An alternative approach would be to use a single branch-and-bound tree, and add the Benders feasibility cuts~\eqref{eq:benders cut} via lazy constraint callbacks. Such an approach is commonly referred to as a ``one-tree" approach, \black{branch-and-cut} {or ``branch-and-check''} and  works much better than its multi-tree counterpart \black{for model \eqref{eq:benders master}}. 

\subsubsection{Branching Strategy}
\label{sec:y form branch}

The branching decisions are extremely important for both the SDP and LP based branch-and-bound algorithms. Following the intuition from Left Branching idea from linear integer programming models as mentioned in Section \ref{sec:lip branch}, we propose a similar scheme that decides the next mark from the left of the ruler. In particular, suppose that the first $m$ marks from the left are located at the positions $\ell_1, \dots, \ell_m$. Then, the location of the mark $m+1$ can be chosen from the set
\[
\left\{ v : \underline{d}_{1,m+1} \le v \le \overline{d}_{1,m+1}, \ v - \ell_k \not\in \{\ell_j - \ell_i: 1 \le i < j \le m\} \ \forall k=1,\dots, m\right\}.
\]
Therefore, we again prefer to create multiple child nodes rather than the more traditional dichotomous branching.

\subsection{Computations}
\label{sec:qip comp}

We report the  results of our computational experiments with the quadratic integer programming formulation with LP linearization in Table~\ref{tab:qipLP} for $n=9,\dots,14$. Since quadratic integer programming formulations are based on feasibility version of the Golomb ruler problem, we solve a sequence of models with increasing ruler length to certify that a given ruler is optimal (see the explanation at the end of Section \ref{sec:feasVers}). For instance, to prove that $G_9=44$, we solve  problem~\eqref{eq:qip feas lp} with $L=35,\dots,43$ and certify that $n_L=8$ (here, we assume  that $G_8$ is known as 34). We report the computational results with and without the Golomb and Cliques cuts giving rise to four different settings.

\begin{table}[t]\footnotesize
\centering
\caption{Results of the ``$y$'' Formulation (linearization via LP) with and without Golomb and Clique cuts using one-tree Benders decomposition. Total time and branch-and-nodes for each mark $n$ and different ruler lengths $L$ are reported.}\label{tab:qipLP}
\begin{tabular}{rccrrcrrcrrcrr}
\toprule
           &            &&         \multicolumn{ 5}{c}{without Golomb Cuts} &&            \multicolumn{ 5}{c}{with Golomb Cuts} \\
\cmidrule{4-8} \cmidrule{10-14}
           &            && \multicolumn{ 2}{c}{w/o Clique Cuts} && \multicolumn{ 2}{c}{with Clique Cuts} && \multicolumn{ 2}{c}{w/o Clique Cuts} && \multicolumn{ 2}{c}{with Clique Cuts} \\
\cmidrule{4-5} \cmidrule{7-8} \cmidrule{10-11} \cmidrule{13-14}
      $   n$ &      $    L$ &&       TT &    \#BBN &&       TT &    \#BBN &&       TT &    \#BBN &&       TT &    \#BBN \\
\midrule
         9 &      35-43 &&       7.91 &      17,753 &&       2.63 &        278 &&       5.73 &      12,630 &&       3.49 &        597 \\

         10 &      45-54 &&      32.24 &      81,036 &&      11.45 &        354 &&      17.73 &      59,770 &&       5.80 &        846 \\

        11 &      56-71 &&    3,850.04 &    1,712,947 &&     730.56 &        824 &&      619.7 &    1,384,860 &&      35.49 &       2,101 \\

        12 &      73-84 &&     $ >10$ hr        &     -       &&    5,576.83 &        674 &&    4,119.82 &    4,710,139 &&      57.80 &       1,661 \\

        13 &     86-105 &&    $ >10$ hr         &     -       &&    $ >10$ hr         &    -        &&      $ >10$ hr       &      -      &&    1,689.52 &       4,869 \\

        14 &    107-126 &&   $ >10$ hr          &     -       &&   $ >10$ hr          &    -        &&     $ >10$ hr        &      -      &&   12,960.17 &       7,776 \\
\bottomrule
\end{tabular}
\end{table}

We observe that both types of cuts are quite effective to solve the subproblems from different perspectives: Golomb cuts  are especially helpful in reducing the computational time more directly whereas Clique cuts significantly lowers the number of branch-and-bound nodes and indirectly reduces the total time. The reason that these two cuts behave differently is that Golomb cuts are added from scratch and their number is limited whereas Clique cuts are added on the fly at each node of the tree and their number can be quite large. We believe that the separation of    Clique cuts can be made more efficient and selective so that the total computational effort can  be further improved.

Finally,  we report the  results of our computational experiments with the quadratic integer programming formulation with SDP linearization in Table~\ref{tab:qipSDP} for $n=9,10,11$ with and without Golomb cuts (MOSEK 8.1 is used as the SDP solver). Although the  total number of branch-and-bound nodes is  reduced by solving the SDP relaxation of the model~\eqref{eq:qip feas sdp} at each node, the total time increases quite significantly which prevents this line of research to be practical.  However, we  point out that our implementation is quite na\"ive and perhaps the value of stronger relaxations provided by the SDP relaxations can be made useful.

\begin{table}[t]\footnotesize
\centering
\caption{Results of the ``$y$'' Formulation (linearization via SDP) with and without Golomb cuts.}\label{tab:qipSDP}
\begin{tabular}{rccrrcrr}
\toprule
           &            &~~& \multicolumn{ 2}{c}{w/o Golomb Cuts} &~~& \multicolumn{ 2}{c}{with Golomb Cuts} \\
\cmidrule{4-5} \cmidrule{7-8}
         $n$ &          $L$ &&       TT &    \#BBN &&       TT &    \#BBN \\
\midrule
         9 &      35-43 &&      21.38 &        227 &&      22.58 &        209 \\

        10 &      45-54 &&      33.03 &        231 &&      40.77 &        200 \\

        11 &      56-71 &&    6,776.11 &      25,591 &&    6,607.07 &      24,737 \\
\bottomrule
\end{tabular}  
\end{table}

\section{Concluding Remarks}

In this paper, we provided a comprehensive \black{comparison} of computational methods to solve  the Golomb ruler problem using optimization techniques. In particular, we analyzed three formulations based on linear integer programming, constraint programming and quadratic integer programming, and proposed several enhancements based on valid inequalities, variable bounding and branching strategies.  According to our experiments with  a budget of 10 hours, integer linear programming models can solve up to 13-mark instances whereas constraint programming and quadratic integer programming formulations can scale up to 14-mark instance, with the latter being faster. We observed that proposed enhancements significantly alter the solution procedures and provide substantial savings in terms of computational effort.

Although the methods in this paper can solve relatively small-size instances of the Golomb ruler problem, we think that there are some promising research directions which might utilize them more effectively. As an example,  if a large number of processors is available, then bound tightening subproblems can be parallelized asynchronously so that they can exchange information whenever a new bound is improved. Since the availability of tight variable bounds seems to accelerate the constraint programming solver, this can potentially enable us to solve larger instances.  Another potential line of research  would be to make the cut generation procedure for the quadratic integer programming model more efficient and selective so that the overhead associated with solving  large MILPs  is reduced while keeping the strength of the relaxations intact.

%
%
%
 \bibliographystyle{splncs04}
 \bibliography{references}
\end{document}